\documentclass[11pt]{amsart}
\usepackage{amsmath,amssymb,amsfonts,latexsym,graphicx,amsthm,mathdots,bbm,enumerate}
\usepackage{hyperref}
\usepackage{url}

\newtheorem{theorem}{Theorem}[section]
\newtheorem{conjecture}[theorem]{Conjecture}
\newtheorem{corollary}[theorem]{Corollary}
\newtheorem{lemma}[theorem]{Lemma}
\newtheorem{proposition}[theorem]{Proposition}
\theoremstyle{definition}
\newtheorem*{remark}{Remark}

\newcommand\bbC{\mathbb{C}}
\newcommand\bbF{\mathbb{F}}
\newcommand\bbk{\mathbbm{k}}
\newcommand\bbQ{\mathbb{Q}}
\newcommand\bbZ{\mathbb{Z}}
\newcommand\calM{\mathcal{M}}
\newcommand\calO{\mathcal{O}}
\newcommand\Ch{\mathrm{Ch}}
\newcommand\one{\mathbbm{1}}
\newcommand\tr{\mathrm{tr}}

\begin{document}

\title[Flag-transitive projective plane]{New necessary conditions for the existence of finite non-Desarguesian flag-transitive projective planes}

\author[Xia]{Binzhou Xia}
\address{School of Mathematics and Statistics\\The University of Melbourne\\Parkville, VIC 3010\\Australia}
\email{binzhoux@unimelb.edu.au}

\begin{abstract}
This paper studies the existence of finite non-Desarguesian flag-transitive projective planes, giving necessary conditions in terms of polynomial equations over finite fields of characteristic $3$. This sheds light on the longstanding conjecture that every finite flag-transitive projective plane is Desarguesian.


\end{abstract}

\maketitle

\section{Introduction}

A \emph{finite projective plane} of order $n$, where $n\geqslant2$ is an integer, is a point-line incidence structure satisfying:
\begin{enumerate}[{\rm(i)}]
\item each line contains exactly $n+1$ points and each point is contained in exactly $n+1$ points;
\item any two distinct lines intersect in exactly one point and any two distinct points are contained in exactly one line.
\end{enumerate}
The incident point-line pairs are a called \emph{flags}. A permutation on the point set preserving the lines and flags is called a \emph{collineation} or \emph{automorphism}. If the collineation group of a finite projective plane acts transitively on the set of flags, then it is said to be \emph{flag-transitive}.

The finite projective planes coordinatized by a finite field are said to be \emph{Desarguesian} since Moufang revealed their equivalence to certain configurational property named in honor of G.~Desargues (see for example~\cite{HP1973}). The celebrated Ostrom-Wagner Theorem~\cite{OW1959} asserts that every finite $2$-transitive projective plane is Desarguesian, which actually classifies the $2$-transitive projective planes.
Note that $2$-transitive finite projective planes are necessarily flag-transitive because two distinct points determine a line.
There is a longstanding conjecture as follows, which, if proved, would strengthen Ostrom-Wagner Theorem.

\begin{conjecture}\label{conj1}
Every finite flag-transitive projective plane is Desarguesian.
\end{conjecture}

Conjecture~\ref{conj1} has been verified by Feit~\cite{Feit1990} for finite flag-transitive projective planes of order up to $14400008$. Moreover, after difficult work in a number of papers, especially that of Kantor~\cite{Kantor1987}, the proof of this conjecture is reduced by the following proposition to a nonexistence problem. However, the conjecture has been ``wide open"~\cite{Gill2016} up until recently. For more on finite flag-transitive projective planes including the result in Proposition~\ref{prop1}, see~\cite{Thas2003,TZ2008}.

\begin{proposition}\label{prop1}
If there exists a finite non-Desarguesian flag-transitive projective plane of order $n$ with $v$ points, then $v=n^2+n+1$ is prime with $n>8$ even, and the set $H_{v,n}$ of nonzero $n$th powers is a $(v,n+1,1)$-difference set in $\bbF_v^+$.
\end{proposition}

A subset $D=\{a_1,\ldots,a_k\}$ of a finite group $G$ is said to be a \emph{$(|G|,k,\lambda)$-difference set} in $G$ or simply a \emph{difference set}, if for each nonidentity $a\in G$ there are exactly $\lambda$ ordered pairs $(a_s,a_t)\in D\times D$ such that $a_sa_t^{-1}=a$. Let $q=m\ell+1$ be a prime power with integers $m$ and $\ell$ such that $1<m<q-1$, and denote
\[
H_{q,m}=\{a^m: a\in\bbF_q^\times\}.
\]
Then $H_{q,m}$ forms a subgroup of $\bbF_q^\times$ of order $\ell$. If this is a $(q,\ell,\lambda)$-difference set in $\bbF_q^+$, then we call it an \emph{$m$th-cyclotomic difference set}. When not specifying the parameters, we will call it a \emph{cyclotomic difference set} for simplicity. For the literature on cyclotomic difference sets, the reader is referred to~\cite{Xia2018} and the references therein.

By virtue of Proposition~\ref{prop1}, the nonexistence of finite non-Desarguesian flag-transitive projective planes can be possibly proved by the nonexistence of cyclotomic difference sets under some conditions. In fact, it is believed that there exists no $m$th-cyclotomic difference set with $m>8$. This is a folklore conjecture, and has been verified by Thas and Zagier~\cite{TZ2008} for prime fields $\bbF_p$ with $p<10^7$. In this paper, we apply $r$-adic Gauss sums and $r$-adic Gauss periods to study the existence of cyclotomic difference sets, and thereby deduce necessary conditions for the existence of finite non-Desarguesian flag-transitive projective plane in terms of polynomial equations over finite fields. Our main theorem is as follows:

\begin{theorem}\label{thm11}
If there exists a finite non-Desarguesian flag-transitive projective plane of order $n$ with $v$ points, then $v=n^2+n+1$ is prime with $n\equiv8\pmod{24}$ and $n>8$, and the following hold:
\begin{enumerate}[{\rm(a)}]
\item there exist a finite field $\bbF$ of characteristic $3$ and $x,x_1,\dots,x_{n-1}\in\bbF$ such that
\begin{equation*}
\begin{cases}
x_{3s}=x^sx_s^3,\quad s=1,\dots,n-1\\
\sum\limits_{\substack{t=1\\t\neq2s}}^{n-1}(-1)^tx_tx_{2s-t}=2x_{2s},\quad s=1,\dots,\frac{n}{2}-1\\
x_sx_{n-s}=(-1)^s,\quad s=1,\dots,\frac{n}{2}\\
x_sx_{\frac{n}{2}+s}=x_{2s}x_{\frac{n}{2}},\quad s=1,\dots,\frac{n}{2}-1\\
x^n=1,
\end{cases}
\end{equation*}
where subscripts of $x$'s are counted modulo $n$;
\item there exist a positive even integer $\sigma$ dividing $n$ and $y_0,y_1,\dots,y_{n-1}\in\bbF_{3^{n/\sigma}}$ such that
\begin{equation*}
\begin{cases}
y_{s+\sigma}=y_s^3,\quad s=0,1,\dots,n-1\\
y_sy_{\frac{n}{2}+s}=-1,\quad s=0,1,\dots,\frac{n}{2}-1\\
\sum\limits_{t=0}^{n-1}y_ty_{s+t}=0,\quad s=0,1,\dots,\frac{n}{2}-1\\
\sum\limits_{t=0}^{n-1}(-1)^ty_ty_{2s-t}=\left(y_s+y_{\frac{n}{2}+s}\right)\sum\limits_{t=0}^{n-1}(-1)^ty_t,\quad s=0,1,\dots,\frac{n}{2}-1,
\end{cases}
\end{equation*}
where subscripts of $y$'s are counted modulo $n$.
\end{enumerate}
\end{theorem}

Theorem~\ref{thm11} will be proved in Section~\ref{sec2} as an application of the main result (Theorem~\ref{thm8}) in Section~\ref{sec1}.
We conjecture that neither of the two systems of polynomial equations in Theorem~\ref{thm11} (even after deleting some equations) has a solution in finite fields of characteristic $3$, which sheds light on Conjecture~\ref{conj1}. Such discussions will be made in Section~\ref{sec3}.

In~\cite{TZ2008}, Thas and Zagier gave equivalent conditions for the existence of cyclotomic difference sets in prime fields by Fermat surfaces and Gauss periods, and applied the latter to determine cyclotomic difference sets $H_{q,m}$ with $q$ prime and $m<10$ (the Gauss period approach was taken in~\cite[Chapter~5]{BEW1998} to determine cyclotomic difference sets $H_{q,m}$ with $q$ prime and $m<10$ or $m=12$). Their approach is different from Lehmer's (using cyclotomic numbers) in~\cite{Lehmer1953}, which enables them to obtain numerical results for prime $q$ up to $10^7$. Note that, for a fixed $m$, there are infinitely many prime numbers $q$ such that $q=m\ell+1$ for some integer $\ell$. In comparison with~\cite{TZ2008}, Theorem~\ref{thm8} cannot lead to numerical results for many values of $q$ without further progress on the polynomial systems therein, since they get too large to solve as $m$ grows. However, Theorem~\ref{thm8} can be used to deduce new nonexistence results for cyclotomic difference sets $H_{q,m}$ with fixed $m$ and arbitrary prime $q$ of the form $m\ell+1$. This is illustrated at the end of Section~\ref{sec1} (see Theorem~\ref{thm12}).

\section{Preliminaries}

For a prime number $r$, let $\left(\frac{\cdot}{r}\right)$ be the Legendre symbol defined by
\[
\left(\frac{n}{r}\right)=
\begin{cases}
1\quad&\text{if $n$ is a square in $\mathbb{F}_r$}\\
-1\quad&\text{if $n$ is a non-square in $\mathbb{F}_r$}
\end{cases}
\]
for all integers $n$ coprime to $r$. Let $\bbZ_r$, $\bbQ_r$ and $\bbC_r$ be the ring of $r$-adic integers, the field of $r$-adic numbers and the $r$-adic completion of the algebraic closure of $\bbQ_r$, respectively. We use $\bbk$ to denote an arbitrary algebraically closed field of characteristic $0$ throughout this section.

\subsection{Gauss sums}

For an integer $n\geqslant2$,
\[
\{x\in\bbk: x^n=1\}
\]
is a cyclic group of order $n$ under multiplication in $\bbk$, and each generator of this group is called a \emph{primitive $n$th root of unity} (\emph{in $\bbk$}). Under pointwise multiplication, the homomorphisms from $\bbF_q^\times$ to $\bbk^\times$ form an abelian group, which is denoted by $\Ch_\bbk(\bbF_q^\times)$. Denote the identity of $\Ch_\bbk(\bbF_q^\times)$ by $\one$. Then $\one$ is the homomorphism from $\bbF_q^\times$ to $\bbk^\times$ sending every element of $\bbF_q^\times$ to $1\in\bbk^\times$. For $\chi\in\Ch_\bbk(\bbF_q^\times)$, extend $\chi$ to a map from $\bbF_q$ to $\bbk$ by setting
\begin{equation*}
\chi(0)=
\begin{cases}
1\quad&\text{if $\chi=\one$}\\
0\quad&\text{if $\chi\neq\one$}
\end{cases}
\end{equation*}
and call this map a (\emph{$\bbk$-valued}) \emph{multiplicative character} on $\bbF_q$. The \emph{order} of a $\bbk$-valued multiplicative character $\chi$ on $\bbF_q$ is defined to be the order of $\chi$ in the group $\Ch_\bbk(\bbF_q^\times)$.

\begin{remark}
For each element $\chi$ of $\Ch_\bbk(\bbF_q^\times)$ and integer $s$, the element $\chi^s$ of $\Ch_\bbk(\bbF_q^\times)$ is also extended to a multiplicative character on $\bbF_q$, and by $\chi^s(\alpha)$ we mean the image of $\alpha\in\bbF_q$ under $\chi^s$ rather than $(\chi(\alpha))^s$. The equality $\chi^s(\alpha)=(\chi(\alpha))^s$ may not hold after $\chi^s$ is extended to a multiplicative character on $\bbF_q$. For example, if $\chi$ has order $s\geqslant2$ then $\chi^s(0)=1\neq0=(\chi(0))^s$.
\end{remark}

For a $\bbk$-valued multiplicative character $\chi$ on $\bbF_q$ and a primitive $p$th root of unity $\zeta$ in $\bbk$, the (\emph{$\bbk$-valued}) \emph{Gauss sum} $G(\chi,\zeta)$ is defined by
$$
G(\chi,\zeta)=\sum_{\alpha\in\bbF_q}\chi(\alpha)\zeta^{\tr(\alpha)},
$$
where $\tr$ is the trace map from $\bbF_q$ to $\bbF_p$. If $r$ is a prime number, then $\bbC_r$-valued Gauss sums are also called \emph{$r$-adic Gauss sums}. Some facts about Gauss sums are listed in Lemma~\ref{lem3} below, whose proof is identical to that of the case $\bbk=\bbC$ in standard textbooks (for example~\cite{BEW1998,Cohen2007I}).

\begin{lemma}\label{lem3}
Let $\bbk$ be an algebraically closed field of characteristic $0$, let $\chi$ be a $\bbk$-valued multiplicative character of order $n$ on $\bbF_q$, and let $\zeta$ be a primitive $p$th root of unity in $\bbk$. Then the following hold:
\begin{enumerate}[{\rm(a)}]
\item for integer $s$ such that $s\not\equiv0\pmod{n}$,
\[
G(\chi^s,\zeta)G(\chi^{-s},\zeta)=\chi^s(-1)q;
\]
\item (Davenport-Hasse product formula) for integer $s$ and positive divisor $d$ of $n$,
\[
\chi^{ds}(d)\prod_{t=0}^{d-1}G(\chi^{s+nt/d},\zeta)=G(\chi^{ds},\zeta)\prod_{t=1}^{d-1}G(\chi^{nt/d},\zeta).
\]
\end{enumerate}
\end{lemma}

\subsection{Discrete Fourier transform}

Fix a primitive $n$th root of unity $\mu$ in $\bbk$ with integer $n\geqslant2$. For a $\bbk$-valued function $X$ on $\mathbb{Z}/n\mathbb{Z}$, the \emph{discrete Fourier transform} (DFT) of $X$, denoted by $\widehat{X}$, is the $\bbk$-valued function on $\mathbb{Z}/n\mathbb{Z}$ defined by
\[
\widehat{X}(s)=\sum\limits_{t=0}^{n-1}\mu^{-st}X(t).
\]
We give in Lemma~\ref{lem4} two well-known formulae for DFT, the convolution formula and the inverse formula (see for example~\cite[page 36]{Terras1999}).

\begin{lemma}\label{lem4}
Let $\bbk$ be an algebraically closed field of characteristic $0$, and let $n\geqslant2$ be an integer. Then the following hold:
\begin{enumerate}[{\rm(a)}]
\item if $W$, $X$ and $Y$ are $\bbk$-valued functions on $\bbZ/n\bbZ$ with
\[
W(s)=\sum\limits_{t=0}^{n-1}X(t)Y(s-t)
\]
for each $s\in\bbZ$, then $\widehat{W}(s)=\widehat{X}(s)\widehat{Y}(s)$ for each $s\in\bbZ/n\bbZ$;
\item if $W$ is a $\bbk$-valued function on $\bbZ/n\bbZ$, then $W(s)=\widehat{\widehat{W}}(-s)/n$ for each $s\in\bbZ/n\bbZ$.
\end{enumerate}
\end{lemma}

\subsection{Gauss periods}

Let $\chi$ be a $\bbk$-valued multiplicative character of order $m$ on $\bbF_q$, let $\zeta$ be a primitive $p$th root of unity in $\bbk$, and let $\mu$ be a primitive $m$th root of unity in $\bbk$. For each $s\in\bbZ$, let
\begin{equation}\label{eq5}
g_s(\chi,\zeta,\mu)=\sum\limits_{t=1}^{m-1}\mu^{-st}G(\chi^t,\zeta)-1.
\end{equation}
We call $g_s(\chi,\zeta,\mu)/m$ a (\emph{$\bbk$-valued}) \emph{Gauss period} or \emph{cyclotomic period} (see~\cite[\S~10.10]{BEW1998}). It is essentially the DFT of Gauss sums. If $r$ is a prime number, then $\bbC_r$-valued Gauss periods are also called \emph{$r$-adic Gauss periods}.

\begin{lemma}\label{lem12}
Suppose that $m$ is even. Then with $g_s=g_s(\chi,\zeta,\mu)$, the following hold:
\begin{enumerate}[{\rm(a)}]
\item for each integer $s$ such that $s\not\equiv m/2\pmod{m}$,
\[
\sum\limits_{t=0}^{m-1}g_tg_{s+t}=m(1-q);
\]
\item if $\chi(2)=\mu^\theta$, then for each integer $s$,
\[
\sum\limits_{t=0}^{m-1}(-1)^tg_tg_{2s-2\theta-t}=\left(g_s+g_{\frac{m}{2}+s}\right)\sum\limits_{t=0}^{m-1}(-1)^tg_t.
\]
\end{enumerate}
\end{lemma}

\begin{proof}
For each $s\in\mathbb{Z}$, let
\[
Y(s)=
\begin{cases}
G(\chi^s,\zeta)\quad&\text{if $s\not\equiv0\pmod{m}$}\\
-1\quad&\text{if $s\equiv0\pmod{m}$}.
\end{cases}
\]
Then $Y$ is a function on $\bbZ/m\bbZ$, and
\[
g_s=\sum\limits_{t=0}^{m-1}\mu^{-st}Y(t)
\]
for all $s$. By Lemma~\ref{lem4}(b) we have
\begin{equation}\label{eq11}
Y\left(\frac{m}{2}\right)=\frac{1}{m}\widehat{\widehat{Y}}\left(-\frac{m}{2}\right)
=\frac{1}{m}\sum\limits_{t=0}^{m-1}\left(\mu^{m/2}\right)^t\widehat{Y}(t)=\frac{1}{m}\sum\limits_{t=0}^{m-1}(-1)^tg_t.
\end{equation}

First suppose that $s$ is an integer with $s\not\equiv m/2\pmod{m}$. Then $-\mu^{-s}\neq1$, and by Lemma~\ref{lem3}(a) we have
\begin{align*}
\sum\limits_{t=0}^{m-1}g_tg_{s+t}
&=\sum\limits_{t=0}^{m-1}\sum\limits_{u=0}^{m-1}\mu^{-tu}Y(u)\sum\limits_{v=0}^{m-1}\mu^{-(s+t)v}Y(v)\\
&=\sum\limits_{u=0}^{m-1}\sum\limits_{v=0}^{m-1}\sum\limits_{t=0}^{m-1}\mu^{-t(u+v)}\mu^{-sv}Y(u)Y(v)\\
&=\sum\limits_{v=0}^{m-1}m\mu^{-sv}Y(-v)Y(v)\\
&=m+\sum\limits_{v=1}^{m-1}m\mu^{-sv}G(\chi^{-v},\zeta)G(\chi^v,\zeta)\\
&=m+\sum\limits_{v=1}^{m-1}m\left(-\mu^{-s}\right)^vq\\
&=m-mq.
\end{align*}
This proves part~(a) of the lemma.

Next suppose that $\chi(2)=\mu^\theta$ and $s$ is an integer. For each integer $v$ such that $1\leqslant v\leqslant m/2-1$ or $m/2+1\leqslant v\leqslant m-1$, we derive from Lemma~\ref{lem3}(b) that
\[
\mu^{2\theta v}G(\chi^v,\zeta)G(\chi^{m/2+v},\zeta)=G(\chi^{2v},\zeta)G(\chi^{m/2},\zeta).
\]
Hence
\[
\mu^{2\theta v}Y\left(\frac{m}{2}+v\right)Y(v)=Y\left(\frac{m}{2}\right)Y(2v)
\]
for each integer $v$. As a consequence,
\begin{align}
\sum\limits_{t=0}^{m-1}(-1)^tg_tg_{2s-2\theta-t}
&=\sum\limits_{t=0}^{m-1}(-1)^t\sum\limits_{u=0}^{m-1}\mu^{-tu}Y(u)\sum\limits_{v=0}^{m-1}\mu^{-(2s-2\theta-t)v}Y(v)\nonumber\\
&=\sum\limits_{u=0}^{m-1}\sum\limits_{v=0}^{m-1}\sum\limits_{t=0}^{m-1}\mu^{t(m/2-u+v)}\mu^{(2\theta-2s)v}Y(u)Y(v)\nonumber\\
&=\sum\limits_{v=0}^{m-1}m\mu^{(2\theta-2s)v}Y\left(\frac{m}{2}+v\right)Y(v)\nonumber\\
&=\sum\limits_{v=0}^{m-1}m\mu^{(2\theta-2s)v}\cdot\mu^{-2\theta v}Y\left(\frac{m}{2}\right)Y(2v)\nonumber\\
&=mY\left(\frac{m}{2}\right)\sum\limits_{v=0}^{m-1}\mu^{-2sv}Y(2v),\nonumber
\end{align}
and so by~\eqref{eq11} we deduce that
\[
\sum\limits_{t=0}^{m-1}(-1)^tg_tg_{2s-2\theta-t}=\sum\limits_{t=0}^{m-1}(-1)^tg_t\sum\limits_{v=0}^{m-1}\mu^{-2sv}Y(2v).
\]
Since
\begin{align*}
\sum\limits_{v=0}^{m-1}\mu^{-2sv}Y(2v)&=2\sum\limits_{\substack{t=0\\ t\text{ even}}}^{m-1}\mu^{-st}Y(t)\\
&=\sum\limits_{t=0}^{m-1}(1+(-1)^t)\mu^{-st}Y(t)\\
&=\sum\limits_{t=0}^{m-1}\mu^{-st}Y(t)+\sum\limits_{t=0}^{m-1}\mu^{-(m/2+s)t}Y(t)\\
&=g_s+g_{\frac{m}{2}+s},
\end{align*}
this leads to
\[
\sum\limits_{t=0}^{m-1}(-1)^tg_tg_{2s-2\theta-t}=\sum\limits_{t=0}^{m-1}(-1)^tg_t\left(g_s+g_{\frac{m}{2}+s}\right),
\]
which proves part~(b) of the lemma.
\end{proof}

\section{Existence of cyclotomic difference sets}\label{sec1}

Throughout this section, let $q=m\ell+1$ be a power of a prime number $p$ with integers $m$ and $\ell$ such that $1<m<q-1$. We study the existence of cyclotomic difference sets, which will yield a system of polynomial equations in Theorem~\ref{thm8} as a necessary condition.

Given a $(v,k,\lambda)$-difference set, we obtain instantly by simple counting that
\[
k(k-1)=\lambda(v-1).
\]
Hence a necessary condition for the existence of cyclotomic $(q,\ell,\lambda)$-difference set is
\begin{equation}\label{eq12}
\ell-1=\lambda m.
\end{equation}

\begin{lemma}\label{lem2}
Let $\bbk$ be an algebraically closed field of characteristic $0$, and let $\chi$ be a $\bbk$-valued multiplicative character of order $m$ on $\bbF_q$. If $m$ is even and $\ell-1=\lambda m$ for some integer $\lambda$, then $\chi(-1)=-1$. In particular, if $H_{q,m}$ is a cyclotomic difference set, then $m$ is even and $\chi(-1)=-1$.
\end{lemma}

\begin{proof}
Suppose that $m$ is even and $\ell-1=\lambda m$ for some integer $\lambda$. Then for any generator $\alpha$ of $\bbF_q^\times$, since $\chi(\alpha)$ is a primitive $m$th root of unity, we have $\chi(\alpha)^{m/2}=-1$. It follows that
\[
\chi(-1)=\chi(\alpha^{(q-1)/2})=\chi(\alpha^{m\ell/2})=(-1)^\ell=(-1)^{\lambda m+1}=-1,
\]
which proves the first conclusion of the lemma.

If $H_{q,m}$ is a $(q,\ell,\lambda)$-difference set, then~\eqref{eq12} holds, and~\cite[Theorem~4.1(a)]{Xia2018} shows that $m$ is even. Thus the second conclusion of the lemma follows.
\end{proof}

In the next lemma we give two sufficient and necessary conditions, in terms of Gauss sums and Gauss periods respectively, for the existence of cyclotomic difference sets. The condition in terms of Gauss sums has appeared in~\cite{Xia2018}. In the special case where the field $\bbF_q$ is prime, other forms of sufficient and necessary conditions in terms of Gauss periods have been established in~\cite[Chapter~5]{BEW1998} and~\cite{TZ2008} to give new proofs of some classical results on cyclotomic difference sets.

\begin{lemma}\label{thm6}
Let $\bbk$ be an algebraically closed field of characteristic $0$, let $\chi$ be a $\bbk$-valued multiplicative character of order $m$ on $\bbF_q$ with $m$ even, let $\zeta$ be a primitive $p$th root of unity in $\bbk$, and let $\mu$ be a primitive $m$th root of unity in $\bbk$. Suppose that $\ell-1=\lambda m$ for some integer $\lambda$. Then the following are equivalent:
\begin{enumerate}[{\rm(i)}]
\item $H_{q,m}$ is a $(q,\ell,\lambda)$-difference set in $\mathbb{F}_q^+$;
\item $\displaystyle\sum_{\substack{t=1\\ t\neq s}}^{m-1}\chi^t(-1)G(\chi^t,\zeta)G(\chi^{s-t},\zeta)=(1+\chi^s(-1))G(\chi^s,\zeta)\,$ for $s=1,\dots,m-1$;
\item $\displaystyle\,g_s(\chi,\zeta,\mu)g_{\frac{m}{2}+s}(\chi,\zeta,\mu)=1+(m-1)q\,$ for each integer $s$.
\end{enumerate}
\end{lemma}

\begin{proof}
The equivalence between~(i) and~(ii) is given in~\cite[Theorem~3.3]{Xia2018}. Thus we only need to prove the equivalence between~(ii) and~(iii). For each $s\in\mathbb{Z}$, let
\[
Y(s)=
\begin{cases}
G(\chi^s,\zeta)\quad&\text{if $s\not\equiv0\pmod{m}$}\\
-1\quad&\text{if $s\equiv0\pmod{m}$}
\end{cases}
\]
and
\[
X(s)=(-1)^sY(s).
\]
Clearly, $X$ and $Y$ are both functions on $\bbZ/m\bbZ$. Define a function on $\bbZ/m\bbZ$ by letting
\[
W(s)=\sum\limits_{t=0}^{m-1}X(t)Y(s-t)
\]
for each $s\in\bbZ$. Then by Lemma~\ref{lem3}(a),
\begin{align*}
W(0)=\sum\limits_{t=0}^{m-1}X(t)Y(-t)&=1+\sum\limits_{t=1}^{m-1}(-1)^tG(\chi^t,\zeta)G(\chi^{-t},\zeta)\\
&=1+\sum\limits_{t=1}^{m-1}q=1+(m-1)q.
\end{align*}
Note $\widehat{Y}(s)=g_s(\chi,\zeta,\mu)$ and $\widehat{X}(s)=g_{\frac{m}{2}+s}(\chi,\zeta,\mu)$. Thus Lemma~\ref{lem4}(a) asserts that
\begin{equation}\label{eq17}
\widehat{W}(s)=\widehat{X}(s)\widehat{Y}(s)=g_s(\chi,\zeta,\mu)g_{\frac{m}{2}+s}(\chi,\zeta,\mu)
\end{equation}
for each integer $s$.

First suppose that~(ii) holds. Then~(i) holds, and so Lemma~\ref{lem2} shows that $m$ is even and $\chi(-1)=-1$. Hence~(ii) turns to
\begin{equation}\label{eq20}
\sum\limits_{t=0}^{m-1}X(t)Y(s-t)=0\quad\text{for}\quad s=1,\dots,m-1.
\end{equation}
Consequently,
\[
W(1)=\dots=W(m-1)=0.
\]
Now for each integer $s$,
\[
\widehat{W}(s)=\sum\limits_{t=0}^{m-1}\mu^{-st}W(t)=W(0)=1+(m-1)q.
\]
This leads to~(iii) by~\eqref{eq17}.

Next suppose that~(iii) holds. It then follows from~\eqref{eq17} that $\widehat{W}(s)=1+(m-1)q$ for each integer $s$. Accordingly, Lemma~\ref{lem4}(b) implies that
\[
W(s)=\frac{1}{m}\sum\limits_{t=0}^{m-1}\mu^{st}\widehat{W}(t)=\frac{1+(m-1)q}{m}\sum\limits_{t=0}^{m-1}\mu^{st}=0\quad\text{for}\quad s=1,\dots,m-1,
\]
which leads to~\eqref{eq20}. Since $m$ is even, we see as in the proof of Lemma~\ref{lem2} that $\chi(-1)=-1$. Thus we obtain~(ii) from~\eqref{eq20}.
\end{proof}

Take a prime number $r$ other than $p$. In what follows we apply Lemma~\ref{thm6} with $\bbk=\bbC_r$. Let $\chi$ be an $r$-adic multiplicative character of order $m$ on $\bbF_q$, let $\zeta$ be a primitive $p$th root of unity in $\bbC_r$, and let $\mu$ be a primitive $m$th root of unity in $\bbC_r$. It is clear from the definition of Gauss sums and Gauss periods that
\[
G(\chi^s,\zeta),g_s(\chi,\zeta,\mu)\in\bbZ[\mu,\zeta]
\]
for each integer $s$. Let $|\ |_r$ be an $r$-adic absolute value on $\bbQ_r(\mu,\zeta)$,
\[
\calO=\{z\in\bbQ_r(\mu,\zeta):|z|_r\leqslant1\}
\]
be the ring of integers of $\bbQ_r(\mu,\zeta)$, and let
\[
\calM=\{z\in\bbQ_r(\mu,\zeta):|z|_r<1\}
\]
be the unique maximal ideal of $\calO$. Then the residue field $\calO/\calM$ is a finite field of characteristic $r$. For any $z\in\calO$, let
\begin{equation}\label{eq4}
\overline{z}=z+\calM\in\calO/\calM.
\end{equation}
Since $\bbZ[\mu,\zeta]\subseteq\calO$, we may consider $\overline{G(\chi^s,\zeta)}$ and $\overline{g_s(\chi,\zeta,\mu)}$ for $s\in\bbZ$.

\begin{lemma}\label{lem11}
Let $r$ be a prime number not dividing $qm$. Then in the above notation, the following hold:
\begin{enumerate}[{\rm(a)}]
\item for each integer $s$ such that $s\not\equiv0\pmod{m}$,
\[
\overline{G(\chi^{rs},\zeta)}=\left(\overline{\chi^r(r)}\right)^s\left(\overline{G(\chi^s,\zeta)}\right)^r;
\]
\item if $\chi(r)=\mu^\sigma$, then for each integer $s$,
\[
\left(\overline{g_s(\chi,\zeta,\mu)}\right)^r=\overline{g_{s+\sigma}(\chi,\zeta,\mu)}.
\]
\end{enumerate}
\end{lemma}

\begin{proof}
Let $s$ be an integer such that $s\not\equiv0\pmod{m}$. Since $\gcd(r,m)=1$, we have $\overline{\chi^{rs}(0)}=0$. Then as $\gcd(r,q)=1$ and $\calO/\calM$ is a field of characteristic $r$, we derive that
\begin{align*}
\overline{G(\chi^{rs},\zeta)}&=\overline{\sum_{\alpha\in\bbF_q^\times}\chi^{rs}(\alpha)\zeta^{\tr(\alpha)}}\\
&=\overline{\sum_{\alpha\in\bbF_q^\times}\chi^{rs}(r\alpha)\zeta^{\tr(r\alpha)}}\\
&=\overline{\sum_{\alpha\in\bbF_q^\times}\chi^{rs}(r)\chi^{rs}(\alpha)\left(\zeta^{\tr(\alpha)}\right)^r}\\
&=\left(\overline{\chi^r(r)}\right)^s\sum_{\alpha\in\bbF_q^\times}\left(\overline{\chi^s(\alpha)\zeta^{\tr(\alpha)}}\right)^r\\
&=\left(\overline{\chi^r(r)}\right)^s\left(\overline{G(\chi^s,\zeta)}\right)^r.
\end{align*}
Hence part~(a) holds. Now suppose $\chi(r)=\mu^\sigma$. It follows that
\[
\left(\overline{G(\chi^t,\zeta)}\right)^r=\left(\overline{\chi(r)}\right)^{-rt}\overline{G(\chi^{rt},\zeta)}
=\overline{\mu}^{\,-\sigma rt}\overline{G(\chi^{rt},\zeta)}\quad\text{for}\quad t=1,\dots,m-1.
\]
Then for each integer $s$ we have
\begin{align*}
\left(\overline{g_s(\chi,\zeta,\mu)}\right)^r&=\left(\sum_{t=1}^{m-1}\overline{\mu^{-st}G(\chi^t,\zeta)}-1\right)^r\\
&=\sum_{t=1}^{m-1}\overline\mu^{\,-rst}\left(\overline{G(\chi^t,\zeta)}\right)^r-1\\
&=\sum_{t=1}^{m-1}\overline\mu^{\,-(s+\sigma)rt}\overline{G(\chi^{rt},\zeta)}-1\\
&=\sum_{t=1}^{m-1}\overline{\mu^{-(s+\sigma)rt}G(\chi^{rt},\zeta)}-1\\
&=\overline{g_{s+\sigma}(\chi,\zeta,\mu)}.
\end{align*}
This proves part~(b).
\end{proof}

To conclude this section, we give the following theorem as necessary conditions from local fields $\bbC_r$ for the existence of cyclotomic difference sets.

\begin{theorem}\label{thm8}
Let $r$ be a prime number not dividing $qm$, let $\chi$ be an $r$-adic multiplicative character of order $m$ on $\bbF_q$, and let $\mu$ be a primitive $m$th root of unity in $\bbC_r$ with $\chi(r)=\mu^\sigma$ and $\chi(2)=\mu^\theta$. Suppose that $H_{q,m}$ is a cyclotomic difference set. Then $m$ is even, and the following hold:
\begin{enumerate}[{\rm(a)}]
\item there exist a finite field $\bbF$ of characteristic $r$ and $u,w,x_1,\dots,x_{m-1}\in\bbF$ such that
\begin{equation*}
\begin{cases}
x_{rs}=u^{\sigma s}x_s^r,\quad s=1,\dots,m-1\\
\sum\limits_{\substack{t=1\\t\neq2s}}^{m-1}(-1)^tx_tx_{2s-t}=2x_{2s},\quad s=1,\dots,\frac{m}{2}-1\\
x_sx_{m-s}=(-1)^sq,\quad s=1,\dots,\frac{m}{2}\\
w^{\theta s}x_sx_{\frac{m}{2}+s}=x_{2s}x_{\frac{m}{2}},\quad s=1,\dots,\frac{m}{2}-1\\
u^m=1\\
w^\frac{m}{2}=1,
\end{cases}
\end{equation*}
where subscripts of $x$'s are counted modulo $m$;
\item there exist $y_0,y_1,\dots,y_{m-1}\in\bbF_{r^{m/\gcd(\sigma,m)}}$ such that
\begin{equation*}
\begin{cases}
y_{s+\sigma}=y_s^r,\quad s=0,1,\dots,m-1\\
y_sy_{\frac{m}{2}+s}=1+(m-1)q,\quad s=0,1,\dots,\frac{m}{2}-1\\
\sum\limits_{t=0}^{m-1}y_ty_{s+t}=m(1-q),\quad s=0,1,\dots,\frac{m}{2}-1\\
\sum\limits_{t=0}^{m-1}(-1)^ty_ty_{2s-2\theta-t}=\left(y_s+y_{\frac{m}{2}+s}\right)\sum\limits_{t=0}^{m-1}(-1)^ty_t,\quad s=0,1,\dots,\frac{m}{2}-1,
\end{cases}
\end{equation*}
where subscripts of $y$'s are counted modulo $m$.
\end{enumerate}
\end{theorem}

\begin{proof}
From Lemma~\ref{lem2} we see that $m$ is even. Let $\zeta$ be a primitive $p$th root of unity in $\bbC_r$. With the notation in~\eqref{eq4}, take $\bbF=\calO/\calM$, $u=\overline{\mu}^{\,r}$, $w=\overline{\mu}^{\,2}$, $x_s=\overline{G(\chi^s,\zeta)}$ for $s=1,\dots,m-1$ and $y_s=\overline{g(\chi^s,\zeta,\mu)}$ for $s=0,1,\dots,m-1$. Then $\bbF$ is a finite field of characteristic $r$.

For the system of equations in part~(a), the first two lines follow from Lemma~\ref{lem11}(a) and Lemma~\ref{thm6}(ii), respectively, and the next two lines follow from Lemma~\ref{lem3}. Since $\mu^m=1$, we have $u^m=(\overline{\mu}^{\,r})^m=\overline{\mu^{rm}}=1$ and $w^{m/2}=(\overline{\mu}^{\,2})^{m/2}=\overline{\mu^m}=1$. This proves part~(a).

For the system of equations in part~(b), the first two lines follow from Lemma~\ref{lem11}(b) and Lemma~\ref{thm6}(iii), respectively, and the last two lines follow from Lemma~\ref{lem12}. Moreover, the first line implies $y_{s+\sigma j}=y_s^{r^j}$ for each nonnegative integer $j$. Hence
\[
y_s=y_{s+\sigma\cdot\frac{m}{\gcd(\sigma,m)}}=y_s^{r^{m/\gcd(\sigma,m)}}
\]
for all $s$, which implies that $y_0,y_1,\dots,y_{m-1}\in\bbF_{r^{m/\gcd(\sigma,m)}}$. This proves part~(b).
\end{proof}

Applying Theorem~\ref{thm8} with a feasible prime number $r$ gives a necessary condition for $m$th-cyclotomic difference sets in terms of polynomial equations over finite fields. We illustrate this by an example of $r=3$. First we present a lemma for general $r$.

\begin{lemma}\label{lem13}
Let $r$ be a prime number not dividing $qm$, let $\chi$ be an $r$-adic multiplicative character of order $m$ on $\bbF_q$, and let $\mu$ be a primitive $m$th root of unity in $\bbC_r$. Suppose that $\chi(r)=\mu^\sigma$ and $m$ is even. Then
\[
(-1)^\sigma=\left(\frac{(-1)^{m/2}q}{r}\right),
\]
where the right-hand side is a Legendre symbol.
\end{lemma}

\begin{proof}
Taking $s=m/2$ in Lemma~\ref{lem3}(a) we obtain
\[
\left(\overline{G(\chi^{m/2},\zeta)}\right)^2=(-1)^{m/2}q
\]
in the notation of~\eqref{eq4}. In particular, $\overline{G(\chi^{m/2},\zeta)}\neq0$. Similarly, by Lemma~\ref{lem11}(a) we have
\[
\left(\overline{G(\chi^{m/2},\zeta)}\right)^r=\left(\overline{\chi(r)}\right)^{-rs}\overline{G(\chi^{rm/2},\zeta)}
=\overline{\mu}^{\,-\sigma rm/2}\overline{G(\chi^{m/2},\zeta)}.
\]
Hence
\[
\left(\overline{G(\chi^{m/2},\zeta)}\right)^{r-1}=\overline{\mu}^{\,-\sigma rm/2}=(-1)^\sigma.
\]
Therefore,
\begin{align*}
(-1)^\sigma=\left(\overline{G(\chi^{m/2},\zeta)}\right)^{r-1}&=\left(\left(\overline{G(\chi^{m/2},\zeta)}\right)^2\right)^{(r-1)/2}\\
&=\left((-1)^{m/2}q\right)^{(r-1)/2}=\left(\frac{(-1)^{m/2}q}{r}\right).\qedhere
\end{align*}
\end{proof}

Now we deduce a consequence of Theorem~\ref{thm8} with $r=3$. This can be used to prove the nonexistence of $m$th-cyclotomic difference sets for some $m$ coprime to $3$.

\begin{corollary}\label{prop2}
Suppose that $H_{q,m}$ is a cyclotomic difference set with $q$ prime, $m>2$ and $m\equiv\pm2\pmod{6}$. Then there exist a positive integer $\sigma$ dividing $m$, an integer $\theta\in\{1,2,\dots,m/2\}$ and $y_0,y_1,\dots,y_{m-1}\in\bbF_{3^{m/\sigma}}$ such that
\begin{equation}\label{eq18}
\begin{cases}
y_{s+\sigma}=y_s^3,\quad s=0,1,\dots,m-1\\
y_sy_{\frac{m}{2}+s}=1+(m-1)(-1)^{\sigma+m/2},\quad s=0,1,\dots,\frac{m}{2}-1\\
\sum\limits_{t=0}^{m-1}y_ty_{s+t}=m(1-(-1)^{\sigma+m/2}),\quad s=0,1,\dots,\frac{m}{2}-1\\
\sum\limits_{t=0}^{m-1}(-1)^ty_ty_{2s-2\theta-t}=\left(y_s+y_{\frac{m}{2}+s}\right)\sum\limits_{t=0}^{m-1}(-1)^ty_t,\quad s=0,1,\dots,\frac{m}{2}-1,
\end{cases}
\end{equation}
where subscripts of $y$'s are counted modulo $m$.
\end{corollary}

\begin{proof}
Let $\chi$ be a $3$-adic multiplicative character of order $m$ on $\bbF_q$, and take $\mu$ to be a primitive $m$th root of unity in $\bbC_3$ such that $\chi(3)=\mu^\sigma$ for some positive integer $\sigma$ dividing $m$. Since $qm$ is not divisible by $3$, it follows from Theorem~\ref{thm8} that there exist an integer $\theta$ and $y_0,y_1,\dots,y_{m-1}\in\bbF_{3^{m/\sigma}}$ with
\begin{equation}\label{eq19}
\begin{cases}
y_{s+\sigma}=y_s^3,\quad s=0,1,\dots,m-1\\
y_sy_{\frac{m}{2}+s}=1+(m-1)q,\quad s=0,1,\dots,\frac{m}{2}-1\\
\sum\limits_{t=0}^{m-1}y_ty_{s+t}=m(1-q),\quad s=0,1,\dots,\frac{m}{2}-1\\
\sum\limits_{t=0}^{m-1}(-1)^ty_ty_{2s-2\theta-t}=\left(y_s+y_{\frac{m}{2}+s}\right)\sum\limits_{t=0}^{m-1}(-1)^ty_t,\quad s=0,1,\dots,\frac{m}{2}-1,
\end{cases}
\end{equation}
where subscripts of $y$'s are counted modulo $m$. Note that~\eqref{eq19} does not change if we replace $\theta$ with $\theta+jm/2$ for any integer $j$. Thus we may assume $\theta\in\{1,2,\dots,m/2\}$.

By Lemma~\ref{lem13},
\[
(-1)^\sigma=\left(\frac{(-1)^{m/2}q}{3}\right).
\]
Accordingly, $(-1)^{m/2}q\equiv(-1)^\sigma\pmod{3}$, which means $q=(-1)^{\sigma+m/2}\in\bbF_{3^{m/\sigma}}$. Substituting this into~\eqref{eq19} we obtain~\eqref{eq18}, as desired.
\end{proof}

Computing the Gr\"{o}bner basis for the system~\eqref{eq18} of polynomial equations in \textsc{Magma}~\cite{BCP1997} we find that, for $m=26$ and $32$ respectively, it has no solution in any field of characteristic $3$ for any positive divisor $\sigma$ of $m$ and $\theta\in\{1,2,\dots,m/2\}$. This yields the following result by Corollary~\ref{prop2}.

\begin{theorem}\label{thm12}
There is neither $26$th nor $32$nd-cyclotomic difference set in any prime field.
\end{theorem}

It would be hopeful that for the family $m\equiv2\pmod{6}$ we are able to obtain more nonexistence results from Corollary~\ref{prop2}, but the author's PC fails to compute the Gr\"{o}bner basis of~\eqref{eq18} for the next value $m=38$ in this family.

\section{Proof of Theorem~\ref{thm11}}\label{sec2}

We have seen in Proposition~\ref{prop1} that the existence problem of finite non-Desarguesian flag-transitive projective plane is related to the existence problem of certain cyclotomic difference sets. Suppose that $H_{v,n}$ is a $(v,n+1,1)$-difference set in $\bbF_v^+$ with $v=n^2+n+1$ prime, as in Proposition~\ref{prop1}. Then clearly, $n$ is not congruent to $1$ modulo $v$, and
\[
n^3=(n^2+n+1)(n-1)+1=v(n-1)+1\equiv1\pmod{v}.
\]
Thereby we conclude that $n$ has order $3$ in $\bbF_v^\times$. Appealing to the First Multiplier Theorem (see~\cite[Theorem~2.1]{Jungnickel1992}) and~\cite[Theorem~IV]{Lehmer1953}) we see that every prime divisor of $n=(n+1)-1$ lies in $H_{v,n}$. This has two consequences. Firstly, $n$ lies in $H_{v,n}$ so that $n^{n+1}=1\in\bbF_v$, which implies that $n+1$ is divisible by $3$. Secondly, noticing that $2$ is a prime divisor of $n$ as $n$ is even, we obtain $2\in H_{v,n}$. To sum up, we have the following:

\begin{lemma}\label{lem10}
If $H_{v,n}$ is a $(v,n+1,1)$-difference set in $\bbF_v^+$ with $v=n^2+n+1$ prime, then $n\equiv2\pmod{3}$ and $2\in H_{v,n}$.
\end{lemma}


We are now in a position to prove Theorem~\ref{thm11}.

\begin{proof}[Proof of Theorem~$\ref{thm11}$]
Suppose that there exists a finite non-Desarguesian flag-transitive projective plane of order $n$ with $v$ points. Then Proposition~\ref{prop1} shows that $v=n^2+n+1$ is prime and $H_{v,n}$ is a $(v,n+1,1)$-difference set in $\bbF_v^+$ with $n>8$ even. By Lemma~\ref{lem10} we have $n\equiv2\pmod{3}$ and $2\in H_{v,n}$. Moreover, we derive from~\cite[Theorem~3.5]{JV1984} that $n\equiv0\pmod{8}$. Hence $n\equiv8\pmod{24}$.

Let $\chi$ be a $3$-adic multiplicative character of order $n$ on $\bbF_v$, and take $\mu$ to be a primitive $n$th root of unity in $\bbC_3$ such that $\chi(3)=\mu^\sigma$ for some positive integer $\sigma$ dividing $n$. Since $2\in H_{v,n}$, it follows that $\chi(2)=1$. Then according to Theorem~\ref{thm8}, there exist a finite field $\bbF$ of characteristic $3$ and $u,x_1,\dots,x_{m-1}\in\bbF$ such that
\begin{equation}\label{eq13}
\begin{cases}
x_{3s}=u^{\sigma s}x_s^3,\quad s=1,\dots,n-1\\
\sum\limits_{t=1}^{2s-1}(-1)^tx_tx_{2s-t}+\sum\limits_{t=2s+1}^{n-1}(-1)^tx_tx_{n+2s-t}=2x_{2s},\quad s=1,\dots,\frac{n}{2}-1\\
x_sx_{n-s}=(-1)^sv,\quad s=1,\dots,\frac{n}{2}\\
x_sx_{\frac{n}{2}+s}=x_{2s}x_{\frac{n}{2}},\quad s=1,\dots,\frac{n}{2}-1\\
u^n=1,
\end{cases}
\end{equation}
and there exist $y_0,y_1,\dots,y_{n-1}\in\bbF_{3^{n/\sigma}}$ such that
\begin{equation}\label{eq14}
\begin{cases}
y_{s+\sigma}=y_s^3,\quad s=0,1,\dots,n-1\\
y_sy_{\frac{n}{2}+s}=1+(n-1)v,\quad s=0,1,\dots,\frac{n}{2}-1\\
\sum\limits_{t=0}^{n-1}y_ty_{s+t}=n(1-v),\quad s=0,1,\dots,\frac{n}{2}-1\\
\sum\limits_{t=0}^{n-1}(-1)^ty_ty_{2s-t}=\left(y_s+y_{\frac{n}{2}+s}\right)\sum\limits_{t=0}^{n-1}(-1)^ty_t,\quad s=0,1,\dots,\frac{n}{2}-1,
\end{cases}
\end{equation}
where subscripts of $y$'s are counted modulo $n$. Since $n\equiv2\pmod{3}$, we have
\[
v=n^2+n+1\equiv2^2+2+1\equiv1\pmod{3}.
\]
Then taking $x=u^\sigma$, the systems~\eqref{eq13} and~\eqref{eq14} give rise to the systems in parts~(a) and~(b) of Theorem~\ref{thm11}, respectively. Finally, as Lemma~\ref{lem13} asserts
\[
(-1)^\sigma=\left(\frac{(-1)^{n/2}v}{3}\right)=\left(\frac{v}{3}\right)=\left(\frac{1}{3}\right)=1,
\]
we deduce that $\sigma$ is even. This completes the proof.
\end{proof}

\section{Discussion on the theorems}\label{sec3}

We make the following two conjectures stating that neither of the two systems of polynomial equations in parts~(a) and~(b) of Theorem~\ref{thm11} has a solution.

\begin{conjecture}\label{conj2}
For any integer $n>8$ with $n\equiv8\pmod{24}$ and any finite field $\bbF$ of characteristic $3$, there do not exist $x,x_1,\dots,x_{n-1}\in\bbF$ such that
\begin{equation}\label{eq26}
\begin{cases}
x_{3s}=x^sx_s^3,\quad s=1,\dots,n-1\\
\sum\limits_{\substack{t=1\\t\neq2s}}^{n-1}(-1)^tx_tx_{2s-t}=2x_{2s},\quad s=1,\dots,\frac{n}{2}-1\\
x_sx_{n-s}=(-1)^s,\quad s=1,\dots,\frac{n}{2}\\
x_sx_{\frac{n}{2}+s}=x_{2s}x_{\frac{n}{2}},\quad s=1,\dots,\frac{n}{2}-1\\
x^n=1,
\end{cases}
\end{equation}
where subscripts of $x$'s are counted modulo $n$.
\end{conjecture}

\begin{conjecture}\label{conj3}
For any integer $n>8$ with $n\equiv8\pmod{24}$ and any positive even integer $\sigma$ dividing $n$, there do not exist $y_0,y_1,\dots,y_{n-1}\in\bbF_{3^{n/\sigma}}$ such that
\begin{equation}\label{eq25}
\begin{cases}
y_{s+\sigma}=y_s^3,\quad s=0,1,\dots,n-1\\
y_sy_{\frac{n}{2}+s}=-1,\quad s=0,1,\dots,\frac{n}{2}-1\\
\sum\limits_{t=0}^{n-1}y_ty_{s+t}=0,\quad s=0,1,\dots,\frac{n}{2}-1\\
\sum\limits_{t=0}^{n-1}(-1)^ty_ty_{2s-t}=\left(y_s+y_{\frac{n}{2}+s}\right)\sum\limits_{t=0}^{n-1}(-1)^ty_t,\quad s=0,1,\dots,\frac{n}{2}-1,
\end{cases}
\end{equation}
where subscripts of $y$'s are counted modulo $n$.
\end{conjecture}

By Theorem~\ref{thm11}, an affirmative answer to either of the above two conjectures will imply the nonexistence of finite non-Desarguesian flag-transitive projective plane and hence confirm Conjecture~\ref{conj1}.

The systems~\eqref{eq26} and~\eqref{eq25} are verified to have no solution when $n=32$. Observe that~\eqref{eq26} has $n$ variables and $(5n-4)/2$ equations, while~\eqref{eq25} has $n$ variables and $5n/2$ equations. As $n$ grows, the systems get more and more over-determined. Hence Conjectures~\ref{conj2} and~\ref{conj3} both seem plausible. In fact, as a stronger version of Conjecture~\ref{conj2}, we conjecture in the following that the system still has no solution if we delete the first line and the last line in~\eqref{eq26}.

\begin{conjecture}
For any integer $n>8$ with $n\equiv8\pmod{24}$ and any finite field $\bbF$ of characteristic $3$, there do not exist $x_1,\dots,x_{n-1}\in\bbF$ such that
\begin{equation*}
\begin{cases}
\sum\limits_{t=1}^{2s-1}(-1)^tx_tx_{2s-t}+\sum\limits_{t=2s+1}^{n-1}(-1)^tx_tx_{n+2s-t}=2x_{2s},\quad s=1,\dots,\frac{n}{2}-1\\
x_sx_{n-s}=(-1)^s,\quad s=1,\dots,\frac{n}{2}\\
x_sx_{\frac{n}{2}+s}=x_{2s}x_{\frac{n}{2}},\quad s=1,\dots,\frac{n}{2}-1.
\end{cases}
\end{equation*}
\end{conjecture}

Similarly, the next conjecture is a stronger version of Conjecture~\ref{conj3} by deleting the first line in~\eqref{eq25}.

\begin{conjecture}\label{conj4}
For any integer $n>8$ with $n\equiv8\pmod{24}$, there do not exist $y_0,y_1,\dots,y_{n-1}\in\bbF_{3^{n/2}}$ such that
\begin{equation*}
\begin{cases}
y_sy_{\frac{n}{2}+s}=-1,\quad s=0,1,\dots,\frac{n}{2}-1\\
\sum\limits_{t=0}^{n-1}y_ty_{s+t}=0,\quad s=0,1,\dots,\frac{n}{2}-1\\
\sum\limits_{t=0}^{n-1}(-1)^ty_ty_{2s-t}=\left(y_s+y_{\frac{n}{2}+s}\right)\sum\limits_{t=0}^{n-1}(-1)^ty_t,\quad s=0,1,\dots,\frac{n}{2}-1,
\end{cases}
\end{equation*}
where subscripts of $y$'s are counted modulo $n$.
\end{conjecture}

We even conjecture that the field $\bbF_{3^{n/2}}$ in Conjecture~\ref{conj4} can be replaced by any field of characteristic $3$.

Let us now turn to results in Section~\ref{sec1}. In the case when $\sigma=1$, the first line of~\eqref{eq18} implies that $y_s=y_0^{r^s}$ for all integer $s\geqslant0$ and thus~\eqref{eq18} is essentially a univariate system on $y_0\in\bbF_{3^m}$. For example, if $\sigma=1$ and $m/2$ is odd, then the first three lines of~\eqref{eq18} gives
\[
\begin{cases}
y_0^{1+3^{m/2}}-m=0\\
\sum\limits_{t=0}^{m-1}y_0^{(1+3^s)3^t}=0,\quad s=0,1,\dots,\frac{m}{2}-1,
\end{cases}
\]
which seems enough to conclude that there is no solution $y_0\in\bbF_{3^m}$. In fact, computation results of the greatest common divisor of the polynomials suggests:

\begin{conjecture}
For an odd prime $r$ and integer $n\geqslant2$, the system of equations
\[
\begin{cases}
y^{1+r^n}-2n=0\\
\sum\limits_{t=0}^{2n-1}y^{(1+r^s)r^t}=0,\quad s=0,1,\dots,n-1
\end{cases}
\]
has a solution $y\in\bbF_{r^{2n}}$ if and only $r$ divides $n$.
\end{conjecture}

Finally, we mention that the ideas in Section~\ref{sec1} might work for topics close to cyclotomic difference sets, such as modified cyclotomic difference sets~\cite{Xia2018}, cyclotomic partial difference sets~\cite{SW2002} and cyclotomic supplementary difference sets~\cite{Yamada1992}, to obtain some results similar to those in Section~\ref{sec1}.

%

\bibliographystyle{plain}

\end{document}